\documentclass[12pt,twoside]{amsart} 
 \title[Abundance theorem]{Abundance theorem for numerically trivial log canonical divisors of semi-log canonical pairs}
\author{Yoshinori  Gongyo}
\address{Graduate School of Mathematical Sciences, 
the University of Tokyo, 3-8-1 Komaba, Meguro-ku, Tokyo 153-8914, Japan.}
\email{gongyo@ms.u-tokyo.ac.jp}
\date{2010/9/11, version 6.02}

\newcommand{\Nklt}[0]{{\operatorname{Nklt}}}

\usepackage{latexsym}
\usepackage{amsmath}
\usepackage{amssymb}
\usepackage{amsthm}
\usepackage{amscd}
\usepackage{enumerate} 
\usepackage{amssymb} 
\usepackage{mathrsfs}          
\usepackage[all]{xy}

\usepackage{latexsym}

\newtheorem{thm}{Theorem}[section]

\newtheorem{prop}[thm]{Proposition}

\newtheorem{lem}[thm]{Lemma}

\newtheorem{rem}[thm]{Remark}

\newtheorem{conj}[thm]{Conjecture}

\newtheorem{cl}[thm]{Claim}
 
\theoremstyle{definition}

\newtheorem{defi}[thm]{Definition}

\newtheorem{step}{Step}

\newtheorem*{ack}{Acknowledgments} 

\newtheorem{theorema}{Theorem}

\begin{document}
\bibliographystyle{amsalpha+}
 
 \maketitle

\begin{abstract}
We prove the abundance theorem for numerically trivial log canonical divisors of log canonical pairs and semi-log canonical pairs.
\end{abstract}

\tableofcontents

\section{Introduction}Throughout this paper, we work over $\mathbb{C}$, the complex number field. We will make use of the standard notation and definitions as in \cite{kmm}. The abundance conjecture is the following:

\begin{conj}[Abundance conjecture]\label{abundance} Let $(X,\Delta)$ be a projective log canonical pair. Then $\nu(K_X+\Delta)=\kappa(K_X+\Delta)$. Moreover, $K_X+\Delta$ is semi-ample if it is nef. 
\end{conj}
For the definition of $\nu(K_X+\Delta)$ and $\kappa(K_X+\Delta)$, we refer \cite{N}. The numerical Kodaira dimension $\nu$ is denoted as $\kappa_{\sigma}$ in \cite{N}. In this paper, we do not use these definitions. The above conjecture is a very important conjecture in the minimal model theory. Indeed, Conjecture \ref{abundance} implies the minimal model conjecture (cf. \cite{b1}, \cite{b2}).
Conjecture \ref{abundance} also says that every minimal model is of general type or has a structure of a Calabi-Yau fiber space, where a Calabi-Yau variety means that its canonical divisor is $\mathbb{Q}$-linearly trivial.
Conjecture \ref{abundance} in dimension $\leq3$ is proved by Fujita, Kawamata, Miyaoka, Keel, Matsuki, and $\mathrm{M^{c}}$Kernan (cf. \cite{k2}, \cite{k1}, \cite{kemamc}). Moreover, Nakayama proved the conjecture when $(X,\Delta)$ is klt and $\nu(K_X+\Delta)=0$ (cf. \cite{N}). Recently, Simpson's result (cf. \cite{simp}) seems to be effective for the proof of the conjecture (cf. \cite{cvp}, \cite{cpt}, \cite{k-ab}, \cite{siu}).   

In this paper, we consider Conjecture \ref{abundance} in the case where $(X,\Delta)$ is minimal and $\nu(K_X+\Delta)=0$, i.e.,  $K_X+\Delta \equiv 0$.

This case for a klt pair is a special case of Nakayama's result. Ambro also gave another proof of the conjecture for klt pairs in this case by using the higher dimensional canonical bundle formula (cf. \cite{A}).
Hence we consider the conjecture in the case where $(X,\Delta)$ is log canonical and $K_X+\Delta \equiv 0$. In Section \ref{lc-ab}, we prove the following theorem by using their results and \cite{BCHM}:

\begin{thm}[=Theorem \ref{main1}]\label{nonvan1}Let $(X,\Delta)$ be a projective log canonical pair. Suppose that $K_X+\Delta \equiv 0$. Then $K_X+\Delta \sim_{\mathbb{Q}}0$.
\end{thm}
  
Moreover, we consider the abundance conjecture for semi-log canonical pairs. 

\begin{defi}[Semi-log canonical]\label{slc}(\cite[Definition 1.1]{F-ab}). Let $X$ be a reduced $S_2$-scheme. We assume that it is pure $n$-dimensional and is normal crossing in codimension $1$. Let $\Delta$ be an effective $\mathbb{Q}$-Weil divisor on $X$ such that $K_X+\Delta$ is $\mathbb{Q}$-Cartier. 

Let $X=\bigcup X_i$ be the decomposition into irreducible components, and $\nu:X':=\coprod X'_i \rightarrow X=\bigcup X_i$ the {\em normalization}, where the \emph{normalization} $\nu : X'=\coprod X'_i \rightarrow X=\bigcup X_i$ means that $\nu|_{X_i'} : X'_i \to X_i$ is the normalization for any $i$. We call $X$ a \emph{normal scheme} if $\nu$ is isomorphic. Define the $\mathbb{Q}$-divisor $\Theta$ on $X'$ by $K_{X'}+\Theta=\nu^*(K_X+\Delta)$ and set $\Theta_i=\Theta|_{X_i'}$. 

We say that $(X,\Delta)$ is \emph{semi-log canonical} (for short, \emph{slc}) if $(X_i',\Theta_i)$ is an lc pair for every $i$. Moreover, we call $(X,\Delta)$ a \emph{semi-divisorial log terminal} (for short, \emph{sdlt}) pair if $X_i$ is normal, that is, $X_i'$ is isomorphic to $X_i$, and $(X_i',\Theta_i)$ is dlt for every $i$. 
 \end{defi}

\begin{conj}\label{semi-abund}Let $(X,\Delta)$ be a projective semi-log canonical pair. Suppose that $K_X+\Delta$ is nef. Then $K_X+\Delta$ is semi-ample.
\end{conj}
Comparing Conjecture \ref{semi-abund} to Conjecture \ref{abundance}, we find that Conjecture \ref{semi-abund} is stated only in the case where $K_X+\Delta$ is nef. In general, it seems that the minimal model program for reducible schemes is difficult. Conjecture \ref{semi-abund} is deeply concerned with the approach to the proof of Conjecture \ref{abundance} (cf. \cite{F-ab}, \cite{Fk-ab}).
 Conjecture \ref{semi-abund} is proved in dimension $\leq 3$ by Kawamata, Abramovich, Fong, Koll\'ar, $\mathrm{M^{c}}$Kernan, and Fujino (cf. \cite{k2}, \cite{AFKM}, \cite{F-ab}). We give an affirmative  answer to Conjecture \ref{semi-abund} in the case where $K_X+\Delta \equiv 0$. 

\begin{thm}\label{nonvan2}Let $(X,\Delta)$ be a projective semi-log canonical pair. Suppose that $K_X+\Delta \equiv 0$. Then $K_X+\Delta \sim_{\mathbb{Q}}0$.
\end{thm}

We prove it along the framework of Fujino in \cite{F-ab}. We take the normalization $\nu:X':=\coprod X'_i \rightarrow X=\bigcup X_i$ and a dlt blow-up on each $X'_i$ (Theorem \ref{dltblowup}). We get $\varphi: (Y,\Gamma) \to (X,\Delta)$ such that $(Y,\Gamma)$ is a (not necessarily connected) dlt pair and  $K_Y+\Gamma=\varphi^*(K_X+\Delta)$. We decompose $Y=\coprod Y'_i$ such that $Y_i \to X_i$ be birational  and put $(K_Y+\Gamma)|_{Y_i}=K_{Y_i}+\Gamma_i$ for every $i$. By Theorem \ref{nonvan1}, it holds that $K_Y+\Gamma \sim_{\mathbb{Q}}0$.  Let $m$ be a sufficiently large and divisible positive integer. Here we need to consider when a section $\{s_i\} \in \bigoplus H^0(Y_i, m(K_{Y_i}+\Gamma_i))$ descends to a section of  $H^0(X,m(K_X+\Delta))$. Fujino introduced the notion of {\em pre-admissible} section and {\em admissible} section for constructing such sections by induction (Definition \ref{admissible}). The pre-admissible sections on $Y$ are descending sections. Moreover Fujino introduced \emph{$B$-pluricanonical representation}
$$\rho_{m}:\mathrm{Bir}(X,\Delta) \to \mathrm{Aut}_{\mathbb{C}}(H^0(X,m(K_X+\Delta)))$$
 for this purpose (cf. Definition \ref{B-bir}, Definition \ref{B-repre}). The basic conjecture of $B$-pluricanonical representation is the following:

\begin{conj}[Finiteness of $B$-pluricanonical representations, {cf. \cite[Conjecture 3.1]{F-ab}}]\label{B-finiteness} Let $(X,\Delta)$ be a projective (not necessarily connected) dlt pair. Suppose that $K_X+\Delta$ is nef. Then $\rho_{m}(\mathrm{Bir}(X,\Delta))$ is finite for a sufficiently large and divisible positive integer $m$.
\end{conj}

This conjecture is proved affirmatively in dimension $\leq 2$ by Fujino (cf. \cite[Theorem 3.3, Theorem 3.4]{F-ab}). To show Theorem \ref{nonvan2}, it  suffices to give an affirmative answer to Conjecture \ref{B-finiteness}  in the case where $K_X+\Delta \equiv 0$. By virtue of Theorem \ref{nonvan1}, it turns out to be that we may assume $K_X+\Delta \sim_{\mathbb{Q}}0$. First, in Section \ref{pre-B-pluri}, we prove Conjecture \ref{B-finiteness} in the case where $(X,\Delta)$ is klt and $K_X+\Delta \sim_{\mathbb{Q}}0$ affirmatively. The proof is  
almost the same as the way of Nakamura--Ueno and Sakai (\cite{nu}, \cite{s}). Next, in Section \ref{semi-log canonical case}, we give an affirmative answer to Conjecture \ref{B-finiteness} under the assumption that $K_X+\Delta \sim_{\mathbb{Q}} 0$ and $\llcorner \Delta\lrcorner \not= 0$ by the way of Fujino (Theorem \ref{finiteness3}). Then we also prove the existence of pre-admissible sections for a (not necessarily connected) dlt pair such that $K_X+\Delta \sim_{\mathbb{Q}} 0$ (Theorem \ref{exis_pread}).
In Section \ref{f-section}, we give some applications of Theorem \ref{nonvan1} and Theorem \ref{nonvan2}. In particular, we remove the assumption of the abundance conjecture from one of the main results in \cite{g} (cf. Theorem \ref{gongyo}). 

After the author submitted this paper to the arXiv, he knew that Theorem \ref{nonvan1} is proved by the same argument as \cite{cvp} in the latest version of \cite{k-ab} and \cite{cvp}. However, in our proof of Theorem \ref{nonvan1}, we do not need Simpson's result.

\begin{ack} The author wishes to express his deep gratitude to his supervisor Professor Hiromichi Takagi for various comments and many suggestions. He also would like to thank Professor Osamu Fujino for fruitful discussions and pointing out my mistakes. He wishes to thank Professor Yoichi Miyaoka for his encouragement. He is also indebted to Doctor Shin-ichi Matsumura for teaching him some knowledge of $L^2$-methods. He is grateful to Professor Yujiro Kawamata for many suggestions. In particular, Professor Kawamata suggests removing the assumption of the abundance conjecture from \cite[Theorem 1.7]{g}. He is supported by the Research Fellowships of the Japan Society for the Promotion of Science for Young Scientists.
\end{ack}

\section {Preliminaries} 
In this section, we introduce some notation and lemmas for the proof of Theorem \ref{nonvan1} and Theorem \ref{nonvan2}. For fixing notation, we start by some basic definitions. The following is the definition of singularities of pairs. Remark that the definitions in \cite{kmm} or \cite{km} are slightly different from ours because the base space is not necessarily connected in our definitions.
\begin{defi}
Let $X$ be a pure $n$-dimensional normal scheme and $\Delta$ a $\mathbb{Q}$-Weil divisor on $X$ such that $K_X+\Delta$ is a $\mathbb{Q}$-Cartier divisor. Let $\varphi:Y\rightarrow X$ be a log resolution of $(X,\Delta)$. We set $$K_Y=\varphi^*(K_X+\Delta)+\sum a_iE_i,$$ where $E_i$ is a prime divisor.
The pair $(X,\Delta)$ is called 
\begin{itemize}
\item[(a)] \emph{sub kawamata log terminal} $($\emph{subklt}, for short$)$ if $a_i > -1$ for all $i$, or
\item[(b)]\emph{sub log canonical} $($\emph{sublc}, for short$)$ if $a_i \geq -1$ for all $i$.
\end{itemize}
If $\Delta$ is effective, we simply call it a \emph{klt} (resp. \emph{lc}) pair. Moreover, we call $X$ a \emph{log terminal} (resp. \emph{log canonical}) variety when $(X,0)$ is \emph{klt} (resp. \emph{lc}) and $X$ is connected.

\end{defi}

\begin{defi}
Let $X$ be a pure $n$-dimensional normal scheme and $\Delta$ an effective $\mathbb{Q}$-Weil divisor on $X$ such that $K_X+\Delta$ is a $\mathbb{Q}$-Cartier divisor. We set an irreducible decomposition $X=\coprod X_i$ and $\Delta_i=\Delta|_{X_i}$. We call that $(X,\Delta)$ is \emph{divisorial log terminal} (for short \emph{dlt}) if $(X_i, \Delta_i)$ is divisorial log terminal for any $i$, where we use the notion of \emph{divisorial log terminal} in \cite{km} for varieties. 
\end{defi}

Ambro and Nakayama prove the abundance theorem for klt pairs whose log canonical divisors are numerically trivial, i.e.,
\begin{thm}[{cf. \cite[Theorem 4.2]{A}, \cite[4.9. Corollary]{N}}]\label{N-nonvan} Let $(X,\Delta)$ be a projective klt pair. Suppose that $K_X+\Delta \equiv 0$. Then $K_X+\Delta \sim_{\mathbb{Q}}0$.
\end{thm}

Next, we introduce a \emph{dlt blow-up}. The following theorem is originally proved by Hacon: 

\begin{thm}[Dlt blow-up, {\cite[Theorem 10.4]{F-fund}, \cite[Theorem 3.1]{kk}}]\label{dltblowup}
Let $X$ be a normal quasi-projective variety and 
$\Delta$ an effective $\mathbb Q$-divisor on $X$ such 
that $K_X+\Delta$ is $\mathbb Q$-Cartier. Suppose that $(X,\Delta)$ is lc.
Then there exists a projective birational 
morphism $\varphi:Y\to X$ from a normal quasi-projective 
variety with the following properties: 
\begin{itemize}
\item[(i)] $Y$ is $\mathbb Q$-factorial, 
\item[(ii)] $a(E, X, \Delta)= -1$ for every  
$\varphi$-exceptional divisor $E$ on $Y$, and
\item[(iii)] for $$
\Gamma=\varphi^{-1}_*\Delta+\sum _{E: {\text{$\varphi$-exceptional}}}E, 
$$ it holds that  $(Y, \Gamma)$ is dlt and $K_Y+\Gamma=\varphi^*(K_X+\Delta)$.
\end{itemize}

\end{thm}

The above theorem is very useful for studying log canonical singularities (cf. \cite{F-index}, \cite{F-fund},  \cite{g}, \cite{kk}).\\
The following elementary lemma is used when we work MMP by Theorem \ref{bchm} (\cite[Corollary 1.3.3]{BCHM}):

\begin{lem}\label{non-pseudo} Let $X$ be an $n$-dimensional normal projective variety such that $n \geq 1$ and $D$ an $\mathbb{R}$-Cartier divisor. Suppose that there exists a nonzero effective $\mathbb{R}$-Cartier divisor $E$ such that $D\equiv -E$. Then $D$ is not pseudo-effective.
\end{lem}

\begin{proof}We take general ample divisors $H_1, \dots, H_{n-1}$. If $D$ is pseudo-effective, then $(D. \bigcap H_i) \geq 0$. But $(E. \bigcap H_i) > 0$. This is a contradiction.
\end{proof}

\begin{rem}\label{rem-non-pseudo} Lemma \ref{non-pseudo} is not always true for such a relative setting as in \cite{kmm} and \cite{BCHM}. For example, let $\pi:X \to U$ be a projective birational morphism. Then every $\mathbb{R}$-Cartier divisor is $\pi$-pseudo-effective. 

\end{rem}

By Birkar--Cascini--Hacon--$\mathrm{M^c}$Kernan, we see the following:

\begin{thm}[{\cite[Corollary 1.3.3]{BCHM}}]\label{bchm} Let $\pi:X \to U$ be a projective morphism of normal quasi-projective varieties and $(X,\Delta)$ a klt pair. Suppose that $K_X+\Delta$ is not $\pi$-pseudo-effective. Then there exists a birational map $\psi : X \dashrightarrow X'$ such that $\psi$ is a composition of $(K_X+\Delta)$-log flips and $(K_X+\Delta)$-divisorial contractions, and $X'$ is a Mori fiber space for $(X,\Delta)$, i.e. there exists an algebraic fiber space $f: X' \to Y'$ such that $\rho(X'/Y')=1$ and $-(K_{X'}+\Delta')$ is $f$-ample, where $\Delta'$ is the strict transform of $\Delta$. 
\end{thm}

\section {Log canonical case}\label{lc-ab} In this section, we prove the follwing:

\begin{thm}\label{main1}Let $(X,\Delta)$ be a projective log canonical pair. Suppose that $K_X+\Delta \equiv 0$. Then $K_X+\Delta \sim_{\mathbb{Q}}0$.
\end{thm}

\begin{proof} We may assume that $\dim X \geq 1$ and $X$ is connected. By taking a dlt blow-up (Theorem \ref{dltblowup}), we also may assume that $(X,\Delta)$ is a $\mathbb{Q}$-factorial dlt pair. By Theorem \ref{N-nonvan}, we may assume that $\llcorner \Delta\lrcorner \not= 0$. We set 
$$S=\epsilon \llcorner \Delta\lrcorner\ \text{and}\ \Gamma = \Delta -S$$
 for some sufficiently small positive number $\epsilon$. Then $(X, \Gamma)$ is klt, $K_X+\Gamma \equiv -S$ is not pseudo-effective by Lemma \ref{non-pseudo}. By Theorem \ref{bchm}, there exist a composition of $(K_X+\Gamma)$-log flips and $(K_X+\Gamma)$-divisorial contractions 
$$\psi: X \dashrightarrow X',$$
 and a Mori fiber space 
$$f':X' \to Y'$$
 for $(X,\Gamma).$ It holds that  $K_{X'}+\Delta'\equiv 0$, where $\Delta'$ is the strict transform of $\Delta$ on $X'$. By the negativity lemma, it suffices to show that $K_{X'}+\Delta'\sim_{\mathbb{Q}} 0$. We put $S'=\psi _*S$ and $\Gamma'= \psi_*\Gamma$. Because it holds $(S'.C)> 0$ for any $f'$-contracting curve $C$, we conclude that $S'\not=0$ and the support of $S'$ dominates $Y'$. Since $K_{X'}+\Delta'\equiv 0$ and $f'$ is a $(K_{X'}+\Gamma')$-extremal contraction, there exists an effective $\mathbb{Q}$-Cartier divisor $D'$ on $Y'$ such that $K_{X'}+\Delta'\sim_{\mathbb{Q}} f'^*D'$ and $D' \equiv 0$ (cf. \cite[Lemma 3-2-5]{kmm}). We remark that $(X',\Delta')$ is not necessarily dlt, but it is a $\mathbb{Q}$-factorial log canonical pair. Hence we can take a dlt blow-up 
$$\varphi:(X'', \Delta'') \to (X',\Delta')$$
 of $(X',\Delta')$. Since the support of $S'$ dominates $Y'$, there exists an lc center $C''$ of $(X'',\Delta'')$ such that $C''$ dominates $Y''$. Then we see that 
$$K_{C''}+\Delta^{''}_{C''} \sim_{\mathbb{Q}} (f'_{C''})^* D',$$
 where $(K_{X''}+\Delta'')|_{C''}=K_{C''}+\Delta''_{C''}$ and $f'_{C''}=f'|_{C''}$. From induction on the dimension of $X$, it holds that $K_{C''}+\Delta^{''}_{C''} \sim_{\mathbb{Q}} 0$. In particular, we conclude that $D'\sim_{\mathbb{Q}}0$. Thus we see that 
$$K_{X'}+\Delta'\sim_{\mathbb{Q}} 0.$$
 We finish the proof of Theorem $\ref{nonvan1}$.
\end{proof}

The above argument does not necessarily hold as it is for a relative setting (cf. Remark \ref{rem-non-pseudo}).

\section{Finiteness of $B$-pluricanonical representations}\label{pre-B-pluri}

Nakamura--Ueno and Deligne proved the following theorem:

\begin{thm}[Finiteness of pluricanonical representations, {\cite[Theorem 14.10]{u}}]\label{nakamura-ueno} Let $X$ be a compact connected Moishezon complex manifold. Then the image of the group homomorphism
$$\rho_{m}:\mathrm{Bim}(X) \to \mathrm{Aut}_{\mathbb{C}}(H^0(X,mK_X))$$
 is finite for any positive integer $m$, where $\mathrm{Bim}(X)$ is the group of bimeromorphic maps from $X$ to itself.
\end{thm}

In this section, we extend Theorem \ref{nakamura-ueno} to klt pairs under the assumption that their log canonical divisors are $\mathbb{Q}$-linear trivial (Theorem \ref{finiteness1}). This is also a generalization of \cite[Proposition 3.1]{F-index} for a sufficiently large and divisible positive integer $m$. The result is used in the proof of Theorem \ref{nonvan2}. Now, we review the notions of $B$-birational maps and $B$-pluricanonical representations introduced by Fujino (cf. \cite{F-ab}).

\begin{defi}[{\cite[Definition 3.1]{F-ab}}]\label{B-bir} Let $(X,\Delta)$ (resp. $(Y,\Gamma)$) be a pair such that $X$ (resp. $Y$) is a normal scheme with a $\mathbb{Q}$-divisor $\Delta$ (resp. $\Gamma$) such that $K_X+\Delta$ (resp. $K_Y+\Gamma$) is $\mathbb{Q}$-Cartier. We say that a proper birational map $f:(X,\Delta)\dashrightarrow (Y,\Gamma)$ is {\em {$B$-birational}} if there exist a common resolution $\alpha:W \to X$ and $\beta: W \to Y$ such that $\alpha^*(K_X+\Delta)=\beta^*(K_Y+\Gamma)$. This means that it holds that $E=F$ when we put $K_W=\alpha^*(K_X +\Delta)+E$ and $K_W=\beta^*(K_Y+\Gamma)+F$. We put $\mathrm{Bir}(X,\Delta)=\{\sigma | \sigma :(X,\Delta) \dashrightarrow (X,\Delta) \text{ is $B$-birational} \}$.
\end{defi}

\begin{defi}[{\cite[Definition 3.2]{F-ab}}]\label{B-repre} Let $X$ be a pure $n$-dimensional normal scheme and $\Delta$ a $\mathbb{Q}$-divisor, and let $m$ be a nonnegative integer such that $m(K_X+\Delta)$ is Cartier. A $B$-birational map $\sigma \in \mathrm{Bir}(X,\Delta)$ defines a linear automorphism of $H^0(X,m(K_X+\Delta))$. Thus we get the group homomorphism 
$$\rho_{m}:\mathrm{Bir}(X,\Delta) \to \mathrm{Aut}_{\mathbb{C}}(H^0(X,m(K_X+\Delta))).$$
The homomorphism $\rho_{m}$ is called a {\em $B$-pluricanonical representation} for $(X,\Delta)$ .  

\end{defi}

Let $X$ be a pure $n$-dimensional normal scheme and $g: X \dashrightarrow X$ a proper birational (or bimeromorphic) map. Set $X= \coprod_{i=1}^{k} X_i$. The map $g$ defines $\sigma \in \mathcal{S}_k$ such that $g|_{X_i}:X_i \dashrightarrow X_{\sigma(i)}$, where $\mathcal{S}_k$ is the symmetric group of degree $k$. Hence $g^{k !}$ induces $g^{k !}|_{X_i} :X_i \dashrightarrow X_{i}$. By Burnside's theorem (\cite[(36.1) Theorem]{cr}), we remark the following:

\begin{rem}\label{b-irr} For the proof of Conjecture \ref{B-finiteness}, we can check that it suffices to show it under the assumption that $X$ is connected. Moreover, Theorem \ref{nakamura-ueno} for a pure dimensional disjoint union of some compact Moishezon complex manifolds holds.
\end{rem}

Now, we show the finiteness of $B$-pluricanonical representations for klt pairs whose log canonical divisors are $\mathbb{Q}$-linearly trivial. Indeed, this result holds for subklt pairs as follows:

\begin{thm}\label{finiteness1}Let $(X,\Delta)$ be a projective subklt pair. Suppose that $K_X+\Delta \sim_{\mathbb{Q}} 0$. Then $\rho_{m}(\mathrm{Bir}(X,\Delta))$ is a finite group for a sufficiently large and divisible positive integer $m$.
\end{thm}

For the proof of Theorem \ref{finiteness1}, the following integrable condition plays an important role:

\begin{defi}\label{integrale-def} Let $X$ be an $n$-dimensional connected complex manifold and   $\omega$ a meromorphic $m$-ple $n$-form. Let $\{U_{\alpha}\}$ be an open covering of $X$ with holomorphic coordinates 
$$(z_{\alpha}^1,z_{\alpha}^2, \cdots, z_{\alpha}^n).$$
We write 

\begin{equation*}\omega|_{U_{\alpha}}=\varphi_{\alpha}(dz_{\alpha}^1\wedge \cdots \wedge dz_{\alpha}^n)^m,
\end{equation*}
where $\varphi_{\alpha}$ is a meromorphic function on $U_{\alpha}$.
We give $(\omega \wedge \bar{\omega})^{1/m}$ by 
\begin{equation*}(\omega \wedge \bar{\omega})^{1/m}|_{U_{\alpha}}=\left( \frac{\sqrt{-1}}{2\pi} \right)^n |\varphi_{\alpha}|^{2/m}dz_{\alpha}^1\wedge d\bar{z}_{\alpha}^1 \cdots \wedge dz_{\alpha}^n \wedge d\bar{z}_{\alpha}^n.
\end{equation*}
We call that a meromorphic $m$-ple $n$-form $\omega$ is {\em $L^{2/m}$-integrable} if $\int_{X} (\omega \wedge \bar{\omega})^{1/m} < \infty$.
\end{defi}

\begin{lem}\label{hol} Let $X$ be a compact connected complex manifold, $D$ a reduced normal crossing divisor on $X$. Set $U =X \setminus D$. If $\omega$ is an $L^2$-integrable meromorphic $n$-form such that $\omega|_{U}$ is holomorphic, then $\omega$ is a holomorphic $n$-form.
\end{lem}
\begin{proof}See \cite[Theorem 2.1]{s} or \cite[Proposition 16]{k3}.
\end{proof}

\begin{lem}\label{integrable} Let $(X,\Delta)$ be a projective subklt pair such that $X$ is a connected smooth variety and $\Delta$ is a simple normal crossing divisor. Let $m$ be a sufficiently large and divisible positive integer, and let $\omega \in H^0(X, \mathcal{O}_{X}(m(K_X+\Delta)))$ be a meromorphic $m$-ple $n$-form. Then $\omega$ is $L^{2/m}$-integrable. 
\end{lem}

\begin{proof} Since $(X,\Delta)$ is subklt, we may write $\Delta=\sum_{i} a_i\Delta_i$, where $\Delta_i$ is a prime divisor and $a_i <1$. We take a sufficiently large and divisible positive integer $m$ such that $1-1/m >a_i$ and $ma_i$ is an integer for any $i$. Thus $\omega$ is a meromorphic $m$-ple $n$-form with at most ($m-1$)-ple pole along $\Delta_i$ for all $i$. By \cite[Theorem 2.1]{s} and holomorphicity of $\omega|_{U}$, $\int_{X} (\omega \wedge \bar{\omega})^{1/m}=\int_{U} (\omega|_{U} \wedge \bar{\omega}|_{U})^{1/m} < \infty$. 
\end{proof}

By Lemma \ref{integrable}, we see the following proposition by almost the same way as \cite[Proposition 1]{nu}, \cite[Proposition 14.4]{u}, and \cite[Lemma 5.1]{s}. 

\begin{prop}\label{algebraic} Let $(X,\Delta)$ be an $n$-dimensional projective subklt pair such that $X$ is smooth and connected, and $\Delta$ is a simple normal crossing divisor. Let $g \in \mathrm{Bir}(X,\Delta)$ be a $B$-birational map, $m$ a sufficiently large and divisible positive integer, and let $\omega \in H^0(X, m(K_X+\Delta))$ be a nonzero meromorphic $m$-ple $n$-form. Suppose that $g^*\omega= \lambda \omega$ for some $\lambda \in \mathbb{C}$. Then there exists a positive integer $N_{m,\omega}$ such that $\lambda^{N_{m,\omega}}=1$ and $N_{m,\omega}$ does not depend on $g$. 
\end{prop}

In the last part of the proof of Proposition \ref{algebraic}, we can avoid the arguments of \cite[Lemma 5.2]{s} (cf. \cite[Proposition 2]{nu}, \cite[Proposition 14.5]{u}) by using Theorem \ref{nakamura-ueno} directly. For the reader's convenience, we include the proof of Proposition \ref{algebraic}.

\begin{proof}[Proof of Proposition \ref{algebraic}]We consider the projective space bundle 
$$\pi:M:=\mathbb{P}_{X}(\mathcal{O}_{X}(-K_X) \bigoplus \mathcal{O}_X) \to X.$$
Set $\Delta=\Delta^+-\Delta^-$, where $\Delta^+$ and $\Delta^-$ are effective, and have no common components. Let $\{U_{\alpha}\}$ be coordinate neighborhoods of $X$ with holomorphic coordinates $(z_{\alpha}^1,z_{\alpha}^2, \cdots, z_{\alpha}^n)$. Since $\omega \in H^0(X,m(K_X+\Delta))$, we can write $\omega$ locally as 
\begin{equation*}
\omega|_{U_{\alpha}}=\frac{\varphi_{\alpha}}{\delta_{\alpha}}(dz_{\alpha}^1\wedge \cdots \wedge dz_{\alpha}^n)^m,
\end{equation*}
where $\varphi_{\alpha}$ and $\delta_{\alpha}$ are holomorphic with no common factors, and $\frac{\varphi_{\alpha}}{\delta_{\alpha}}$ has poles at most $m\Delta^{+}$. We may assume that $\{U_{\alpha}\}$ gives a local trivialization of $M$, i.e. $M|_{U_{\alpha}} := \pi^{-1}U_{\alpha} \simeq U_{\alpha} \times \mathbb{P}^1$. We set a coordinate $(z_{\alpha}^1,z_{\alpha}^2, \cdots, z_{\alpha}^n, \xi_{\alpha}^1:\xi_{\alpha}^2)$ of $U_{\alpha} \times \mathbb{P}^1$ with homogeneous coordinates $(\xi_{\alpha}^1:\xi_{\alpha}^2)$ of $\mathbb{P}^1$. Note that 
$$\frac{\xi_{\alpha}^1}{\xi_{\alpha}^2}= k_{\alpha \beta} \frac{\xi_{\beta}^1}{\xi_{\beta}^2}\ \text{in}\ M|_{U_{\alpha}\bigcap U_{\beta}},$$
 where $k_{\alpha \beta}=\mathrm{det}(\partial z_{\beta}^i/\partial z_{\alpha}^j)_{1\leq i,j \leq n}.$ Set 
$$Y_{U_{\alpha}}=\{ (\xi_{\alpha}^1)^m\delta_{\alpha}-(\xi_{\alpha}^2)^m\varphi_{\alpha}=0\} \subset U_{\alpha} \times \mathbb{P}^1.$$
We can patch $\{ Y_{U_{\alpha}}\}$ easily and denotes the patching by $Y$. Remark that $Y$ may have singularities and be reducible.  Let $\pi_1:M' \to M$ be a log resolution of $(M, Y\cup \pi^{-1}( \mathrm{Supp} \Delta))$ such that $Y'$ is smooth, where $Y'$ is the strict transform of $Y$ on $M'$. We set $F'= \pi \circ \pi_1$ and $f'= F'|_{Y'}$. Remark that $Y'$ may  be disconnected and a general fiber of $f'$ is $m$ points. Define a meromorphic $n$-form on $M$ by 
$$\Theta|_{M|_{U_{\alpha}}}=\frac{\xi_{\alpha}^1}{\xi_{\alpha}^2}dz_{\alpha}^1\wedge \cdots \wedge dz_{\alpha}^n.$$
We put $\theta'=\pi_1^*\Theta|_{Y'}$. By the definition, 
$$ (\theta')^m=f'^* \omega.$$
Since $\int_{X} (\omega \wedge \bar{\omega})^{1/m} < \infty$ by Lemma \ref{integrable}, it holds that $\int_{Y'} \theta' \wedge \bar{\theta'} < \infty$. Hence $\theta'$ is $L^2$-integrable. Since $f'^{-1}(\mathrm{Supp}\Delta)$ is simple normal crossings, $\theta'$ is a holomorphic $n$-form on $Y'$ by Lemma \ref{hol}. 

We take a $\nu \in \mathbb{R}$ such that $\nu^m=\lambda$. We define a birational map $\bar{g}_{\nu}: M \dashrightarrow M$ by 
$$\bar{g}_{\nu}: (z_{\alpha}^1,z_{\alpha}^2, \cdots, z_{\alpha}^n, \xi_{\alpha}^1:\xi_{\alpha}^2) \to (g(z_{\alpha}^1,z_{\alpha}^2, \cdots, z_{\alpha}^n), \nu (\mathrm{det}(\partial g / \partial z_{\alpha} ))^{-1}\xi_{\alpha}^1: \xi_{\alpha}^2)$$
 on $U_{\alpha}$. Then $\bar{g}_{\nu}$ induces a birational map $h':Y' \dashrightarrow Y'$. It satisfies that

\begin{equation*}
\label{eq:line}
\xymatrix{
Y' \ar[d]_{f'} \ar@{-->}[r]^{h'} \ar@{}[dr]|\circlearrowleft & Y' \ar[d]^{f'} \\
X \ar@{-->}[r]^{g} & X. \\
}\
\end{equation*} 

Thus we see 

$$ h'^*(\theta')^m=h'^*f'^* \omega = f'^* g^* \omega =\lambda f'^* \omega=\lambda(\theta')^m.
$$

Because Theorem \ref{nakamura-ueno} holds for pure dimensional possibly disconnected projective manifolds (Remark \ref{b-irr}), there exists a positive integer $N_{m,\omega}$ such that $\lambda^{N_{m,\omega}}=1$ and $N_{m,\omega}$ does not depend on $g$. We finish the proof of Proposition \ref{algebraic}.

\end{proof}

\begin{proof}[Proof of Theorem \ref{finiteness1}]By taking a log resolution of $(X,\Delta)$, we may assume that $X$ is smooth and $\Delta$ has a simple normal crossing support. Let $m$ be a sufficiently  large and divisible positive integer. Since $\dim_{\mathbb{C}}H^0(X,m(K_X + \Delta))=1$ by the assumption that $K_X+\Delta \sim_{\mathbb{Q}} 0$, we see that $\rho_{m}(g) \in \mathbb{C}^*$ for any $g \in \mathrm{Bir}(X,\Delta)$. Proposition \ref{algebraic} implies that $(\rho_{m}(g))^{N_{m,\omega}}=1$. Hence $\rho_{m}(\mathrm{Bir}(X,\Delta))$ is a finite group because it is a subgroup of $\mathbb{C}^*$.
\end{proof}

\section{Semi-log canonical case}\label{semi-log canonical case}

In this section, following the framework of \cite{F-ab} in the case where $K_X+\Delta \equiv 0$, we prove Theorem \ref{nonvan2}. Here we recall the definition of {\em sdlt} pairs as  in Definition \ref{slc}.

\begin{defi}[{cf.~\cite[Definition 4.1]{F-ab}}]\label{admissible}
Let $(X,\Delta)$ be an $n$-dimensional proper sdlt pair
and $m$ a sufficiently divisible integer. We take the normalization $\nu:X':=\coprod X'_i \rightarrow X=\bigcup X_i$.
We define {\em{admissible}} and {\em{pre-admissible}} sections inductively on dimension as the follows:
\begin{itemize}
\item $s\in H^{0}(X,m(K_X+\Delta))$ is {\em{pre-admissible}} if the restriction 
$$\nu^*s|_{(\amalg _{i}{\llcorner\Theta_{i}\lrcorner})}\in 
H^{0}(\amalg _{i} {\llcorner\Theta_{i}\lrcorner},m(K_{X'}+\Theta)|_
{(\amalg _{i}{\llcorner\Theta_{i}\lrcorner})})$$
is admissible.
\item $s\in H^{0}(X,m(K_X+\Delta))$ is {\em{admissible}}
if $s$ is pre-admissible and $g^{*}(s|_{X_j})=s|_{X_i}$ for 
every $B$-birational map $g:(X_i,\Theta_i)\dashrightarrow (X_j,\Theta_j)$ for 
every $i,j$.
\end{itemize}
Note that if $s\in H^{0}(X,m(K_X+\Delta))$ 
is admissible, then the restriction $s|_{X_i}$ is 
$\mathrm{Bir}(X_i,\Theta_i)$-invariant for every $i$.
\end{defi}

\begin{rem}\label{rem-ad} Let $(X,\Delta)$ be an $n$-dimensional proper sdlt pair
and $m$ a positive integer such that $m(K_X+\Delta)$ is Cartier. We take the normalization $\nu:X' \rightarrow X$. Then it is clear by definition that $s\in H^{0}(X,m(K_X+\Delta))$ is admissible (resp. pre-admissible) if and only if so is $\nu^*s\in H^{0}(X',m(K_{X'}+\Delta'))$.
\end{rem}

For the normalization $\nu:X' \rightarrow X$, any pre-admissible section on $X'$ descends on $X$ (cf. \cite[Lemma 4.2]{F-ab}). Therefore in our case we sufficiently prove the existence of nonzero pre-admissible sections on $X'$. Including this statement, we prove the following three theorems by induction on the dimension:

\begin{theorema}\label{exis_pread} Let $(X,\Delta)$ be an $n$-dimensional projective (not necessarily connected) dlt pair. Suppose that $K_X+\Delta\sim_{\mathbb{Q}}0$. Then there exists a nonzero pre-admissible section $s\in H^{0}(X,m(K_X+\Delta))$ for a sufficiently large and divisible positive integer $m$.
\end{theorema}

\begin{theorema}\label{finiteness3} Let $(X,\Delta)$ be an $n$-dimensional projective (not necessarily connected) dlt pair. Suppose that $K_X+\Delta \sim_{\mathbb{Q}} 0$. Then $\rho_{m}(\mathrm{Bir}(X,\Delta))$ is a finite group for a sufficiently large and divisible positive integer $m$.
\end{theorema}

\begin{theorema}\label{exis_ad} Let $(X,\Delta)$ be an $n$-dimensional projective (not necessarily connected) dlt pair. Suppose that $K_X+\Delta\sim_{\mathbb{Q}}0$. Then there exists a nonzero admissible section $s\in H^{0}(X,m(K_X+\Delta))$ for a sufficiently large and divisible positive integer $m$.

\end{theorema}

\begin{step}\label{induction2} Theorem \ref{exis_ad}$_{n-1}$ implies Theorem \ref{exis_pread}$_{n}$.
\end{step}

we claim the following  by using Theorem \ref{bchm}: 

\begin{cl}[{cf. \cite[12.3.2. Proposition]{AFKM}, \cite[Proposition 2.1]{F-ab}, \cite[Proposition 2.4]{F-index}, \cite[Proposition 5.1]{kk}}]\label{fibration} Let $(X,\Delta)$ be an $n$-dimensional $\mathbb{Q}$-factorial connected dlt pair such that $n \geq 2$. Suppose that $K_X+\Delta \sim_{\mathbb{Q}} 0$. Then one of the following holds: 

\begin{itemize}
\item[(i)] $\llcorner \Delta\lrcorner$ is connected, or
\item[(ii)] $\llcorner \Delta\lrcorner$ has two connected components $\Delta_{1}$ and $\Delta_{2}$. Moreover, there exist a birational map $\varphi:X \dashrightarrow X'$ and an algebraic fiber space 
$f':X' \to  Y'$ with a general fiber $\mathbb{P}^1$  such that they satisfy the following:

\begin{itemize}
\item[(ii-a)] $\varphi$ is a composition of  $(K_X+\Delta)$-log flops and $(K_X+\Delta)$-crepant divisorial contractions, and $X'$ is log terminal,
\item[(ii-b)]$Y'$ is an $(n-1)$-dimensional $\mathbb{Q}$-factorial projective log terminal variety, and
\item[(ii-c)] there exists an effective $\mathbb{Q}$-divisor $\Omega'$ on $Y'$ such that $(Y', \Omega')$ is an lc pair and $f'^*(K_{Y'}+\Omega')=K_{X'}+\Delta'$, where $\Delta'=\varphi_{*}\Delta$.
\end{itemize}

Furthermore, there exists an irreducible component $D_i \subset \Delta_i$ such that $f'|_{D_i'}:(D_i', \Delta'_{D_i'}) \rightarrow (Y',\Omega')$ is a $B$-birational isomorphism for $i=1,2$, where $D_i':= \varphi_{*}D_i$ and $K_{D_i'}+\Delta'_{D_i'} = (K_{X'}+\Delta')|_{D_i'}$. In particular, $(f' \circ \varphi)|_{D_i}:(D_i, \Delta_{D_i}) \dashrightarrow (Y',\Omega')$ is a $B$-birational map, where $K_{D_i}+\Delta_{D_i} = (K_{X}+\Delta)|_{D_i}$.

\end{itemize}
\end{cl}

\begin{proof}
We set $S=\epsilon \llcorner \Delta\lrcorner$ and $\Gamma:=\Delta -S$ for some sufficiently small positive number $\epsilon$. For the proof of Theorem \ref{nonvan1}, we can get a Mori fiber space $f':X' \to Y'$ for $(X,\Gamma)$ such that a birational map $\varphi:X \dashrightarrow X'$ is a composition of $(K_X+\Gamma)$-log flips and $(K_X+\Gamma)$-divisorial contractions. It holds that  $K_{X'}+\Delta'\sim_{\mathbb{Q}} 0$, where $\Delta'$ is the strict transform of $\Delta$ on $X'$. By the proof of \cite[Proposition 2.1]{F-ab}, we see Claim \ref{fibration}.

\end{proof}

Thus the following claim holds for an $n$-dimensional dlt pair $(X,\Delta)$ such that $K_X+\Delta \sim_{\mathbb{Q}}0$ from the same way as \cite[Proposition 4.5]{F-ab} by Claim \ref{fibration}.  

\begin{cl}[{cf. \cite[Proposition 4.5]{F-ab}, \cite[Proposition  4.15]{F-index}}]\label{fujino} 
Let $(X,\Delta)$ be an $n$-dimensional $\mathbb{Q}$-factorial (not necessarily connected) dlt pair. Suppose that $K_X+\Delta \sim_{\mathbb{Q}} 0$. Let $m$ be a sufficiently large and divisible positive integer, and put $S=\llcorner \Delta \lrcorner$. For a nonzero admissible section $s \in H^0(X,m(K_{S}+\Delta_{S}))$, there exists a nonzero pre-admissible section $t \in H^0(X, m(K_X+\Delta))$ such that $s|_{S}=t$, where $K_{S}+\Delta_{S} = (K_X+\Delta)|_{S}$.  In particular, Theorem \ref{exis_ad}$_{n-1}$ implies Theorem \ref{exis_pread}$_{n}$.
\end{cl} 

We finish the proof of Step \ref{induction2}.

Next, we see the following:

\begin{step}\label{induction3} Theorem \ref{exis_pread}$_{n}$ implies Theorem \ref{finiteness3}$_n$.
\end{step}

\begin{proof} (cf. \cite[Theorem 3.5]{F-ab}). By Remark \ref{b-irr}, we may assume that $X$ is connected. We may assume that $\llcorner \Delta\lrcorner \not =0$ by Theorem \ref{finiteness1}. Because we assume that Theorem \ref{exis_pread}$_{n}$ holds, we can take a nowhere vanishing section $s\in H^{0}(X,m(K_X+\Delta))$ for a sufficiently large and divisible positive integer $m$. Since $\dim_{\mathbb{C}}H^0(X,m(K_X + \Delta))=1$, we see that $\rho_{m}(g) \in \mathbb{C}^*$ for any $g \in \mathrm{Bir}(X,\Delta)$. By \cite[Lemma 4.9]{F-ab}, it holds that $\rho_{m}(g)s|_{\llcorner \Delta\lrcorner}=g^*s|_{\llcorner \Delta\lrcorner}=s|_{\llcorner \Delta\lrcorner}$ for any $g \in \mathrm{Bir}(X,\Delta)$. Thus, it holds that $\rho_{m}(g) =1$ for any $g \in \mathrm{Bir}(X,\Delta)$. Hence the action of $\mathrm{Bir}(X,\Delta)$ on $H^{0}(X,m(K_X+\Delta))$ is trivial.
\end{proof}

Lastly, the following step follows directly by \cite[Lemma 4.9]{F-ab}.  

\begin{step}\label{induction1} Theorem \ref{exis_pread}$_{n}$ and Theorem \ref{finiteness3}$_n$ imply Theorem \ref{exis_ad}$_{n}$.
\end{step}

Thus we obtain Theorem \ref{exis_pread}, Theorem \ref{finiteness3}, and Theorem  \ref{exis_ad}.

Finally, we show Theorem \ref{nonvan2}:

\begin{proof}[Proof of Theorem \ref{nonvan2}] We take the normalization $\nu:X':=\coprod X'_i \rightarrow X=\bigcup X_i$ and a dlt blow-up on each $X'_i$. We get $\varphi: (Y,\Gamma) \to (X,\Delta)$ such that $(Y,\Gamma)$ is a (not necessarily connected) dlt pair and  $K_Y+\Gamma=\varphi^*(K_X+\Delta)$. By Theorem \ref{nonvan1}, it holds that $K_Y+\Gamma \sim_{\mathbb{Q}}0$. Thus there exists a nonzero pre-admissible section $s \in H^0(Y,m(K_Y+\Gamma))$ for a sufficiently large and divisible positive integer $m$ by Theorem \ref{exis_pread}. Therefore $s$ descends on $X$ by \cite[Lemma 4.2]{F-ab}. Because the descending section of $s$ is nonzero, it holds that $K_X+\Delta \sim_{\mathbb{Q}}0$. We finish the proof of Theorem \ref{nonvan2}.
\end{proof}

\section{Applications}\label{f-section}

In this section, we give some applications on Theorem \ref{nonvan1} and Theorem \ref{nonvan2}.

First, we expand the following Fukuda's theorem to $4$-dimensional log canonical pairs.

\begin{thm}[{\cite[Theorem\ 0.1]{Fk}}]\label{fk} Let $(X,\Delta)$ be a projective klt pair. Suppose that there exists a semi-ample $\mathbb{Q}$-divisor $D$ such that $K_{X}+\Delta \equiv D$. Then $K_X+\Delta$ is semi-ample. 
\end{thm}

\begin{thm}\label{fukuda} Let $(X,\Delta)$ be a $4$-dimensional projective log canonical pair. Suppose that there exists a semi-ample $\mathbb{Q}$-divisor $D$ such that $K_{X}+\Delta \equiv D$. Then $K_X+\Delta$ is semi-ample. 
\end{thm}

We obtain Theorem \ref{fukuda} from the same argument as \cite{Fk} by replacing Kawamata's theorem (\cite[Theorem 4.3]{k2}, \cite[Theorem 6-1-11]{kmm} \cite[Theorem1.1]{F-kawamata}) with the following theorem:   

\begin{thm}[{\cite[Corollary 6.7]{F-bp}}]\label{fuji-bp}
Let $(X, B)$ be an lc pair and 
$\pi:X \to S$ a proper 
morphism onto a variety $S$. Assume the following 
conditions: 
\begin{itemize}
\item[(a)] $H$ is a $\pi$-nef $\mathbb Q$-Cartier divisor on $X$, 
\item[(b)] $H-(K_X+B)$ is $\pi$-nef and $\pi$-abundant, 
\item[(c)] $\kappa (X_\eta, (aH-(K_X+B))_\eta)\geq 0$ and $\nu(X_\eta, (aH-(K_X+B))_\eta)=\nu(X_\eta, (H-(K_X+B))_\eta)$ for 
some $a\in \mathbb Q$ with $a>1$, 
where $\eta$ is the generic point of $S$,  and
\item[(d)] there is a positive integer $c$ such that 
$cH$ is Cartier and 
$\mathcal O_T(cH|_T)$
 is 
$\pi$-generated, 
where $T=\Nklt(X, B)$ is the non-klt locus of $(X, B)$.  
\end{itemize}
Then $H$ is $\pi$-semi-ample. 
\end{thm}

\begin{proof}[Proof of Theorem \ref{fukuda}] (cf.\ \cite{Fk}). By taking a dlt blow up, we may assume that $(X,\Delta)$ is a $\mathbb{Q}$-factorial dlt pair. Since $D$ is semi-ample, We get an algebraic fiber space $f:X \to Y$ such that $D=f^*A$ for some $\mathbb{Q}$-ample divisor on $Y$. \\
By Theorem \ref{nonvan1} and the abundance theorem for semi-log canonical threefolds (cf. \cite{F-ab}), $f: X \to Y$ satisfies the condition of the assumption in Theorem \ref{fuji-bp}. Thus $K_X+\Delta$ is $f$-semi-ample. We get an algebraic fiber space $g:X \to Z$ over $Y$ such that there exists some $\mathbb{Q}$-ample divisor $B$ over $Y$ such that $K_{X}+\Delta = g^*B$.

Because $K_X+\Delta \equiv D$, it holds that $f\simeq g$ as algebraic fiber spaces.

Since $A \equiv B$, $B$ is ample. Thus we see that $K_X+\Delta$ is semi-ample. We finish the proof of Theorem \ref{fukuda}.
 \end{proof}

Next, as an application of Theorem \ref{nonvan1}, we obtain the following theorem by \cite[Theorem 1.7]{g}: 

\begin{thm}[{\cite[Theorem 1.7]{g}}]\label{gongyo} Let $(X, \Delta)$ be an $n$-dimensional lc weak log Fano pair, that is $-(K_X+\Delta)$ is nef and big, and $(X, \Delta)$ is lc. Suppose that $\mathrm{dim}\mathrm{Nklt}(X, \Delta)\leq 1$. Then $-(K_X+\Delta)$ is semi-ample.
\end{thm}

By the same way as the proof of \cite[Theorem 1.7]{g}, we also see the following:

\begin{thm}\label{gongyo2} Let $(X, \Delta)$ be an $n$-dimensional lc pair such that $K_X+\Delta$ is nef and big. Suppose that $\mathrm{dim}\mathrm{Nklt}(X, \Delta)\leq 1$.
Then $K_X+\Delta$ is semi-ample.
\end{thm}


\begin{thebibliography}{BCHM}

\bibitem[AFKM]{AFKM}
D. Abramovich, L. L. Y. Fong, J. Koll\'ar and J. $\mathrm{M^{c}}$Kernan, Semi log canonical surface, Flip and Abundance for algebraic threefolds, Ast\'erisque 211 (1992), 139--154.

\bibitem[A]{A}
F. Ambro, The moduli $b$-divisor of an lc-trivial fibration,  Compositio. Math.  141  (2005),  no. 2, 385--403.

\bibitem[B1]{b1}

C. Birkar, On existence of log minimal models, preprint, arXiv:0706.1792, to appear in Compositio Math.

\bibitem[B2]{b2}

\bysame, On existence of log minimal models II, preprint, arXiv:0907.4170, to appear in J. Reine Angew Math.

\bibitem[BCHM]{BCHM}
C. Birkar, P. Cascini, C. D. Hacon and J. $\mathrm{M^{c}}$Kernan, Existence of minimal models for varieties of log general type, J. Amer. Math. Soc. 23 (2010), 405-468.
\bibitem[CKP]{cvp}
F.Campana, V. Koziarz and M. P$\mathrm{\breve{a}}$un, Numerical character of the effectivity of adjoint line bundles, preprint, arXiv:1004.0584.

\bibitem[CPT]{cpt}
F. Campana, T. Peternell and M. Toma, Geometric stability of the cotangent bundle and the universal cover of a projective manifold, preprint, math/0405093.

\bibitem[CR]{cr}
C. W Curtis and I.  Reiner, \textit{Representation theory of finite groups and associative algebras}, Reprint of the 1962 original. Wiley Classics Library. A Wiley-Interscience Publication. John Wiley \& Sons, Inc., New York, 1988.

\bibitem[Fj1]{F-ab}
O. Fujino, Abundance theorem for semi log canonical threefolds, Duke Math. J. 102  (2000),  no. 3, 513--532.

\bibitem[Fj2]{F-index}
\bysame, The indices of log canonical singularities, Amer. J. Math.  123  (2001),  no. 2, 229--253. 

\bibitem[Fj3]{F-bp}
\bysame, Base point free theorems--saturation, b-divisors, and canonical bundle formula--, math.AG/0508554. 

\bibitem[Fj4]{F-kawamata}
\bysame, On Kawamata's theorem, preprint, arXiv:0910.1156, to appear in the Schiermonnikoog volume.

\bibitem[Fj5]{F-fund}
\bysame, Fundamental theorems for the log minimal model program, preprint, arXiv:0909.4445.

\bibitem[Fk1]{Fk-ab}
S. Fukuda, On numerically effective log canonical divisors, Int. J. Math. Math. Sci. 30 (2002), no. 9, 521?531.

\bibitem[Fk2]{Fk}
\bysame, An elementary semi-ampleness result for log canonical divisors, preprint, arXiv:1003.1388.

\bibitem[G]{g}
Y. Gongyo, On weak Fano varieties with log canonical singularities, preprint, arXiv:0911.0974.

\bibitem[Ka1]{k3}
Y. Kawamata, Characterization of abelian varieties, Compositio Math.  43  (1981), no. 2, 253--276.

\bibitem[Ka2]{k2}
\bysame, Pluricanonical systems on minimal algebraic varieties, Inv. Math. 79 (1985), no. 3, 567-588.

\bibitem[Ka3]{k1}

\bysame, Abundance theorem for minimal threefolds.  Invent. Math.  108  (1992),  no. 2, 229--246.

\bibitem[Ka4]{k-ab}
\bysame, On the abundance theorem in the case of $\nu=0$, preprint, arXiv:1002.2682

\bibitem[KaMM]{kmm}
Y. Kawamata, K, Matsuda and K, Matsuki, \textit{Introduction to the minimal model problem}, Algebraic geometry, Sendai, 1985,  283--360, Adv. Stud. Pure Math., 10, North-Holland, Amsterdam, 1987. 

\bibitem[KeMM]{kemamc}
S. Keel, K. Matsuki and J. $\mathrm{M^{c}}$Kernan, Log abundance theorem for threefolds.  Duke Math. J.  75  (1994),  no. 1, 99--119, Corrections to: ``Log abundance theorem for threefolds'',  Duke Math. J.  122  (2004),  no. 3, 625--630. 

\bibitem[KoKo]{kk}
J. Koll\'ar and S. J Kov\'acs,  Log canonical singularities are Du Bois, preprint, arXiv:0902.0648.

\bibitem[KoM]{km}
J. Koll\'ar and S. Mori. \textit{Birational geometry of algebraic varieties}, Cambridge Tracts in Math.,134 (1998).

\bibitem[NU]{nu}
I. Nakamura and K. Ueno, An addition formula for Kodaira dimensions of analytic fibre bundles whose fibre are Moi$\mathrm{\check{s}}$ezon manifolds,  J. Math. Soc. Japan  25  (1973), 363--371.

\bibitem[N]{N}
N. Nakayama, \textit{Zarisiki decomposition and abundance}, MSJ Memoirs, 14. Mathematical Society of Japan, Tokyo, 2004. 

\bibitem[S]{s}
F. Sakai, \textit{Kodaira dimensions of complements of divisors,  Complex analysis and algebraic geometry},  pp. 239--257. Iwanami Shoten, Tokyo, 1977.

\bibitem[Sim]{simp}
C. Simpson, Subspaces of moduli spaces of rank one local systems,  Ann. Sci. \'Ecole Norm. Sup. (4)  26  (1993),  no. 3, 361--401. 

\bibitem[Siu]{siu}
Y. T. Siu, Abundance conjecture, preprint, arXiv:0912.0576.

\bibitem[U]{u}
K. Ueno, \textit{Classification theory of algebraic varieties and compact complex spaces}, Lecture Notes in Math., Vol. 439, Springer, Berlin, 1975.


\end{thebibliography}
\end{document}